\documentclass[titlepage,11pt]{article}
\usepackage{latexsym}
\usepackage{amsfonts}
\usepackage{amsmath}
\newcommand\blackslug{\hbox{\hskip 1pt \vrule width 4pt height 8pt depth 1.5pt
        \hskip 1pt}}
\newcommand\bbox{\hfill \quad \blackslug \bigbreak}

\def\lx{,\ldots,}

\title{Induced subgraphs of graphs with large chromatic number. \\
IX. Rainbow paths}
\author{Alex Scott\\
Oxford University, Oxford, UK
\\
\\
Paul Seymour\thanks{Supported by ONR grant N00014-14-1-0084 and NSF
grant DMS-1265563.}\\
Princeton University, Princeton, NJ 08544, USA}

\date{January 20, 2017; revised July 3, 2017}

\newtheorem{thm}{}[section]

\newcommand{\Proof}{\noindent{\bf Proof.}\ \ }

\begin{document}
\maketitle
\begin{abstract}
We prove that for all integers $\kappa, s\ge 0$ there exists $c$ with the following property. Let $G$ be a graph with 
clique number at most $\kappa$ and chromatic number more than $c$. Then for every 
vertex-colouring (not necessarily optimal) of $G$, 
some induced subgraph of $G$ is an $s$-vertex path, and all its vertices have different colours.
This extends a recent result of Gy\'arf\'as and S\'ark\"ozy~\cite{gyarfas2}, who proved the same 
for graphs $G$ with  $\kappa=2$ and girth at least five.

\end{abstract}

\section{Introduction}

Graphs in this paper are finite and have no loops or multiple edges.
We denote the chromatic number and the clique number of $G$ by $\chi(G), \omega(G)$ respectively.
If $X\subseteq V(G)$, the subgraph of $G$ induced on $X$ is denoted by $G[X]$.
A {\em colouring} of a graph $G$ is a map $\phi$ from $V(G)$ to the set of positive integers such that $\phi(u)\ne \phi(v)$
for all adjacent $u,v$; and a {\em coloured graph} is a pair $(G,\phi)$ where $G$ is a graph and $\phi$ is a colouring of $G$.
Given a coloured graph $(G,\phi)$, a subgraph $H$ of $G$ is said to be {\em rainbow} 
if $\phi(u)\ne\phi(v)$ for all distinct $u,v\in V(H)$.  

The following interesting conjecture on rainbow paths is due to Aravind (see~\cite{babu}).
\begin{thm}\label{aravind} 
{\bf Conjecture:} 
Let $G$ be a triangle-free graph.  Then for every colouring (not necessarily optimal) of $G$, there is a rainbow induced subgraph
isomorphic to a $\chi(G)$-vertex path.
\end{thm}
This remains open, but some special cases have been proved. For instance, if we just ask for an induced path
(not necessarily rainbow), then
it holds by a theorem of Gy\'arf\'as~\cite{gyarfas}.  Or if we just ask for a rainbow path 
(not necessarily induced), then it holds by the Gallai-Roy theorem~\cite{gallai, roy}, 
even without the bound on clique number: if we direct every edge of $G$ 
towards the end with higher colour, then every directed path of the digraph obtained is rainbow.
The conjecture also holds if the girth of $G$ equals its chromatic number, by a result of Babu, Basavaraju, Chandran and 
Francis~\cite{babu}: in particular, if $(G,\phi)$ is a triangle-free coloured graph with $\chi(G)\ge 4$, 
then some induced four-vertex path of $G$ is rainbow.

A recent paper of 
Gy\'arf\'as and S\'ark\"ozy~\cite{gyarfas2} proves the following result.
\begin{thm}\label{gs} 
For all $s\ge 1$ there exists $c$ such that the following holds.
Let $G$ be a graph with girth at least 5 and $\chi(G)>c$.  Then for every colouring of $G$ there is a rainbow induced subgraph
isomorphic to an $s$-vertex path.
\end{thm}

In this paper, we extend this theorem in two ways: we remove the girth restriction, and allow a general bound on clique size.
Here is our result.

\begin{thm}\label{mainthm}
For all $\kappa,s \ge 1$ there exists $c$ such that for every coloured graph $(G,\phi)$ with $\omega(G)\le \kappa$
and $\chi(G)>c$, there is a rainbow induced subgraph of $G$ isomorphic to an $s$-vertex path.
\end{thm} 

We prove this in the next section, and include some further discussion in the conclusion.

\section{The proof}

We will need the following theorem of Galvin, Rival and Sands~\cite{galvin}:
\begin{thm}\label{galvin}
For all integers $s\ge 0$ there exists $r\ge 0$ with the following property. For every graph $G$ that has a path with 
at least $r$ vertices, either some induced path of $G$ has at least $s$ vertices, or some subgraph of $G$ is isomorphic
to the complete bipartite graph $K_{s,s}$.
\end{thm}

A {\em grading} in a graph $G$ is a sequence $(W_1\lx W_n)$ of subsets of $V(G)$, pairwise
disjoint and with union $V(G)$. If $w\ge 0$ is such that $\chi(G[W_i])\le w$ for $1\le i\le n$ we say the grading 
is {\em $w$-colourable}. We say that $u\in V(G)$ is {\em earlier} than $v\in V(G)$, and $v$ is {\em later} than $u$
(with respect to some grading $(W_1\lx W_n)$) if $u\in W_i$ and $v\in W_j$ where $i<j$.
We need the following lemma:

\begin{thm}\label{gradinglemma}
Let $s\ge 0$ be an integer, and let $r$ be as in \ref{galvin}. Let $w\ge 0$, and let 
$(G,\phi)$ be a coloured graph with $\chi(G)\ge wr$, such that no $s$-vertex induced path of $G$ is rainbow. 
Let $(W_1\lx W_n)$ be a $w$-colourable grading in $G$. 
Then there exist $i\in \{1\lx n\}$, and a vertex $v\in W_i$, and a set of $s$ vertices, pairwise
with different colours, all later than $v$ and all adjacent to $v$.
\end{thm}
\Proof
Choose $r$ as in \ref{galvin}. Since each $G[W_i]$ is $w$-colourable, there is a partition $(A_1\lx A_w)$ of $V(G)$
such that $W_i\cap A_j$ is stable for $1\le i\le n$ and $1\le j\le w$. 
Since $\chi(G)\ge wr$,  there exists $j$ such that $\chi(G[A_j])>r$.
For each edge $e=uv$ of $G[A_j]$ direct $e$ 
from $u$ to $v$ if $\phi(u)<\phi(v)$,
obtaining a digraph $D$ say. 
By the Gallai-Roy theorem, there is a directed path $P$ of $D$
with $\chi(G[A_j])>r$ vertices. From the definition of $D$, it follows that all vertices of $P$ have different colours.
Since no $s$-vertex induced path of $G$ is rainbow, it follows from \ref{galvin} applied to 
$G[V(P)]$ that some subgraph $H$ of $G[V(P)]$ is isomorphic to $K_{s,s}$ and rainbow. Choose $i\in \{1\lx n\}$
minimum such that $W_i\cap V(H)\ne \emptyset$, and choose $v\in W_i\cap V(H)$. Since $H$ is rainbow, there
are $s$ vertices of $H$ all with different colours and all adjacent to $v$. From the choice of $i$, none of them
is earlier than $v$; and since they all belong to $V(H)\subseteq A_j$ and $W_i\cap A_j$ is stable, 
none of them belongs to $W_j$. Consequently they are all later than $v$. This proves \ref{gradinglemma}.~\bbox

Now we prove \ref{mainthm}, which we restate:

\begin{thm}\label{mainthm2}
Let $s,\kappa\ge 0$ be integers. Then there exists $c\ge 0$ such that for every coloured graph $(G,\phi)$ with 
$\omega(G)\le \kappa$ and $\chi(G)>c$, some induced $s$-vertex path of $G$ is rainbow.
\end{thm}
\Proof
Since the result holds if $\kappa\le 1$, we may assume by induction on $\kappa$ that $\kappa\ge 2$ 
and there exists $c'$ such that for every coloured graph $(G,\phi)$ with
$\omega(G)\le \kappa-1$ and $\chi(G)>c'$, some induced $s$-vertex path of $G$ is rainbow.
Let $r$ be as in \ref{galvin}. 
Define $w_s=0$, and for $j=s-1\lx 0$ let $w_j=w_{j+1}r+c'$. Let $c=(w_1+1)r$; we claim that $c$ satisfies the theorem.
Let $(G,\phi)$ be a coloured graph with $\omega(G)\le \kappa$ and $\chi(G)>c$. We must show
that some induced $s$-vertex path of $G$ is rainbow. Suppose not.
For each vertex $v$, if $N$ denotes the set of neighbours of $v$,
then $\omega(G[N])\le \kappa-1$, and so $\chi(G[N])\le c'$.

For each vertex $z$, let $A(z)$ be 
the set of all vertices $v$ such that there is an induced
rainbow path of $G$ between $z,v$. 
\\
\\
(1) {\em $\chi(G[A(z)])>w_1+1$ for some vertex $z$.}
\\
\\
For suppose not. Let $V(G)=\{v_1\lx v_n\}$, and for $1\le i\le n$ let 
$$W_i=A(v_i)\setminus (A(v_1)\cup \cdots\cup A(v_{i-1})).$$
Thus $(W_1\lx W_n)$ is a $(w_1+1)$-colourable grading in $G$. By \ref{gradinglemma}, since $\chi(G)>c= (w_1+1)r$,
there exist $i\in \{1\lx n\}$, and a vertex $v\in W_i$, and a set $X$ of $s$ vertices, pairwise
with different colours, all later than $v$ and all adjacent to $v$. Since $v\in A(v_i)$, there is an induced rainbow
path of $G$ between $v_i,v$, say $P$. Since $G$ has no induced rainbow $s$-vertex path, $|V(P)|<s$. Consequently
some vertex $x\in X$ has a colour different from the colours of the vertices of $P$. But then adding $x$ to $P$
gives a rainbow path between $v_i,x$, and therefore there is an induced rainbow path between $v_i,x$. Consequently
$x\in A(v_i)$. But $x$ is later than $v$, a contradiction. This proves (1).

\bigskip

Choose $z$ as in (1). Let $Q$ be a rainbow induced path of $G$ with first vertex $z$. An {\em extension}
of $Q$ is a rainbow induced path $P$ of $G$ with first vertex $z$ such that $Q$ is a subpath of $P$ and $Q\ne P$.
We denote by $B(Q)$
the set of all vertices $v\in V(G)$ such that there is an extension of $Q$ between $z,v$.
Choose a rainbow induced path $Q$ of $G$ with first vertex $z$, and $\chi(G[B(Q)])>w_j$
where $j=|V(Q)|$, such that $j$ is maximum. (This is well-defined, since every such path has fewer than $s$ vertices
by hypothesis, and since $Q$ exists with $j=1$.) Thus $j<s$.
Let $Q$ have ends $z,y$ say, and let $V$ be 
the set of all vertices of $G$ adjacent to $y$, adjacent to no other vertex of $Q$, and with a different colour
from every vertex of $Q$. Thus every extension of $Q$ contains a vertex in $V$; and
for each $v\in V$, adding $v$ to $Q$ gives an extension of $Q$,
say $Q_v$. Every vertex in $B(Q)$ belongs either to $V$ or to $B(Q_v)$ for some $v\in V$; and the vertices in $B(Q_v)$
for some $v$ are precisely the vertices in $B(Q)$ that are nonadjacent to $y$, and these vertices have no neighbours
in $V(Q)$ at all. Now 
$\omega(G[V])\le \kappa-1$, and so $\chi(G[V])\le c'$; and consequently
$\chi(B(Q)\setminus V)>w_j-c'=w_{j+1}r$.

From the choice of $Q$, $\chi(G[B(Q_v)])\le w_{j+1}$ for each $v\in V$. Let $v=\{v_1\lx v_n\}$,
and for $1\le i\le n$ let $B_i=B(Q_{v_i})$, and $W_i = B_i\setminus (B_1\cup\cdots\cup B_{i-1})$.
Thus $(W_1\lx W_n)$ is a $w_{j+1}$-colourable grading of $G[B(Q)\setminus V]$. 
By \ref{gradinglemma}, there exist $i\in \{1\lx n\}$ and $v\in W_i$, and a set $X$ of neighbours of $v$,
all later than $v$ and all with different colours, with $|X|=s$. Since $v\in W_i\subseteq B_i$, there is an extension
$P$  of $Q$ between $z,v$, and therefore $|V(P)|<s$ by hypothesis.
Since $|X|\ge s$, there exists $x\in X$ with a colour different from the colour of every vertex in $P$. 
Adding $x$
to $P$ gives a rainbow path between $z,x$ of which $P$ is a subpath; and since $P$ is induced, and $x$ has no neighbour
in $V(Q)$, there is an extension of $Q(v_i)$ between $z,x$, and so 
$x\in B_i$. But this is impossible since $x$ is later than $v$. This contradiction shows that there is
a rainbow induced path in $G$ with $s$ vertices, and so proves \ref{mainthm2}.~\bbox

\section{Conclusion}

Can \ref{mainthm} be extended beyond paths?
One could ask which graphs $H$ have the following property:
for every triangle-free graph $G$ with sufficiently large chromatic number, and for every colouring of $G$,
some induced subgraph of $G$
is isomorphic to $H$ and rainbow. But for $H$ to have this property:
\begin{itemize}
\item $H$ must have no cycles. Otherwise we can take
$G$ to have girth larger than the length of such a cycle and large chromatic number, and then $G$ contains no copy
of $H$ at all, rainbow or otherwise.
\item Every vertex of $H$ must have degree at most two.  Otherwise,
as shown by Kierstead and Trotter~\cite{kierstead}, a counterexample
is given by the ``shift graph'' of triples: the vertex set is $[n]^{(3)}$, and two triples are adjacent 
if the smallest two elements of one are the same as the largest two elements of the other.
If we colour every triple by its middle element then no rainbow subgraph has a vertex of degree more than two. 
\end{itemize}
So the only graphs $H$
which might have the desired property are forests with maximum degree at most two, that is, induced subgraphs of paths.

We might also ask: is it true that if $G$ is triangle-free and has sufficiently large chromatic number then for every colouring of $G$,
some induced subgraph of $G$ is isomorphic to a cycle and is rainbow?  In other words, does $G$ contain a rainbow {\em hole}?  
But here again the shift graph of triples (again coloured by middle elements) 
gives a counterexample.  Indeed, this coloured graph does not even have a hole in which every three consecutive vertices are rainbow
(every three-vertex path is monotonic, so there is no way for the cycle to ``close up'').  But we do not know the following:
is it true that for all fixed $s,\kappa$, if $G$ is a graph with $\omega(G)\le \kappa$ and $\chi(G)$ sufficiently large
then in every colouring of $G$ there is a hole in which {\em some} set of $s$ consecutive vertices is rainbow?
(We note that, resolving an old conjecture of Gy\'arf\'as~\cite{gyarfas}, it was shown in~\cite{css} that every such graph 
does at least contain a long hole.)

Let us also mention another question to which we do not know the answer. Let $G$ be a triangle-free graph with very large
chromatic number, and let $\mathcal{A}$ be a set of stable subsets (not necessarily
pairwise disjoint) of $V(G)$ with union $V(G)$.
Does there necessarily exist an $s$-vertex induced path $P$ of $G$ such that for each
$v\in V(P)$, some $X\in \mathcal{A}$ satisfies $X\cap V(P)=\{v\}$?

Finally, we remark that we do not believe \ref{aravind}, but have not found a counterexample.  Since \ref{aravind}
is known to hold for $\chi(G)=4$, the first place to look for counterexamples to \ref{aravind} is when $\chi(G)= 5$ and we want a 
five-vertex induced rainbow path. In a laborious and unavailing search for a 
counterexample, we checked by hand all colourings of $G$
when $G$ is the Mycielski graph on 23 vertices, for which $\chi(G)= 5$ and 
$\omega(G)=2$; they all satisfy the conjecture.

\end{document}